# 0–1 LAWS FOR REGULAR CONDITIONAL DISTRIBUTIONS

By Patrizia Berti and Pietro Rigo

*Universitá di Modena e Reggio Emilia and Universitá di Pavia*

Let $(\Omega, \mathcal{B}, P)$ be a probability space, $\mathcal{A} \subset \mathcal{B}$ a sub-$\sigma$-field, and $\mu$ a regular conditional distribution for $P$ given $\mathcal{A}$. Necessary and sufficient conditions for $\mu(\omega)(A)$ to be 0–1, for all $A \in \mathcal{A}$ and $\omega \in A_0$, where $A_0 \in \mathcal{A}$ and $P(A_0) = 1$, are given. Such conditions apply, in particular, when $\mathcal{A}$ is a tail sub-$\sigma$-field. Let $H(\omega)$ denote the $\mathcal{A}$-atom including the point $\omega \in \Omega$. Necessary and sufficient conditions for $\mu(\omega)(H(\omega))$ to be 0–1, for all $\omega \in A_0$, are also given. If $(\Omega, \mathcal{B})$ is a standard space, the latter 0–1 law is true for various classically interesting sub-$\sigma$-fields $\mathcal{A}$, including tail, symmetric, invariant, as well as some sub-$\sigma$-fields connected with continuous time processes.

**1. Introduction and motivations.** Let $(\Omega, \mathcal{B}, P)$ be a probability space and $\mathcal{A} \subset \mathcal{B}$ a sub-$\sigma$-field. A *regular conditional distribution* (r.c.d.) for $P$ given $\mathcal{A}$ is a mapping $\mu : \Omega \to \mathbb{P}$, where $\mathbb{P}$ denotes the set of probability measures on $\mathcal{B}$, such that $\mu(\cdot)(B)$ is a version of $E(I_B \mid \mathcal{A})$ for all $B \in \mathcal{B}$. A $\sigma$-field is *countably generated* (c.g.) in case it is generated by one of its countable subclasses. In the sequel, it is assumed that $P$ admits a r.c.d. given $\mathcal{A}$ and

(1)      $\mu$ denotes a *fixed* r.c.d., for $P$ given $\mathcal{A}$, and $\mathcal{B}$ is c.g.

Moreover,

$$H(\omega) = \bigcap_{\omega \in A \in \mathcal{A}} A$$

is the atom of $\mathcal{A}$ including the point $\omega \in \Omega$.

Heuristically, conditioning to $\mathcal{A}$ should mean conditioning to the atom of $\mathcal{A}$ which actually occurs, say $H(\omega)$, and the probability of $H(\omega)$ given $H(\omega)$ should be 1. If this interpretation is agreed, $\mu$ should be everywhere *proper*, that is, $\mu(\omega)(A) = I_A(\omega)$ for all $A \in \mathcal{A}$ and $\omega \in \Omega$. Though $\mu(\cdot)(A) = I_A(\cdot)$ a.s. for *fixed* $A \in \mathcal{A}$, however, $\mu$ can behave quite inconsistently with properness.

Say that $\mathcal{A}$ is c.g. under $P$ in case the trace $\sigma$-field $\mathcal{A} \cap C = \{A \cap C : A \in \mathcal{A}\}$ is c.g. for some $C \in \mathcal{A}$ with $P(C) = 1$. If $\mathcal{A}$ is c.g. under $P$ then

(2) $\quad\quad\quad \mu(\omega)(A) = I_A(\omega) \quad \text{for all } A \in \mathcal{A} \text{ and } \omega \in A_0,$

where, here and in what follows, $A_0$ designates some set of $\mathcal{A}$ with $P(A_0) = 1$. In general, the exceptional set $A_0^c$ cannot be removed. Further, (2) implies that $\mathcal{A} \cap A_0$ is c.g. See [2, 3] and Theorem 1 of [4]. In other terms, not only everywhere properness is to be weakened into condition (2), but the latter holds if and only if $\mathcal{A}$ is c.g. under $P$.

If $\mathcal{A}$ fails to be c.g. under $P$, various weaker notions of r.c.d. have been investigated. Roughly speaking, in such notions, $\mu$ is asked to be everywhere proper but $\sigma$-additivity and/or measurability are relaxed. See [1] and references therein. However, not very much is known on r.c.d.'s, regarded in the usual sense, when $\mathcal{A}$ is not c.g. under $P$ (one exception is [11]). In particular, when (2) fails, one question is whether some of its consequences are still available.

In this paper, among these consequences, we focus on:

(3) $\quad\quad\quad \mu(\omega)(H(\omega)) \in \{0, 1\} \quad \text{for all } \omega \in A_0,$

(4) $\quad\quad\quad \mu(\omega)(A) \in \{0, 1\} \quad \text{for all } A \in \mathcal{A} \text{ and } \omega \in A_0.$

Note that (4) implies (3) in case $H(\omega) \in \mathcal{A}$ for all $\omega \in A_0$, and that, for (3) to make sense, one needs to assume $H(\omega) \in \mathcal{B}$ for all $\omega \in A_0$.

Both conditions (3) and (4) worth some attention.

Investigating (3) can be seen as a development of the seminal work of [4]. The conjecture that (3) holds (under mild conditions) is supported by those examples in the literature where $\mathcal{B}$ is c.g. In these examples, in fact, either $\mu(\omega)(H(\omega)) = 1$ a.s. or $\mu(\omega)(H(\omega)) = 0$ a.s. See, for instance, [11].

Condition (4) seems to have been neglected so far, though it is implicit in some ideas of [7] and [5]. In any case, (4) holds in a number of real situations and can be attached a clear heuristic meaning. As to the latter, fix $\omega_0 \in \Omega$. Since $\mu(\omega_0)$ comes out by conditioning on $\mathcal{A}$, one could expect that $\mu$ is a r.c.d. for $\mu(\omega_0)$ given $\mathcal{A}$, too. Condition (4) grants that this is true, provided $\omega_0 \in A_0$. More precisely, letting $M = \{Q \in \mathbb{P} : \mu \text{ is a r.c.d. for } Q \text{ given } \mathcal{A}\}$, condition (4) is equivalent to

$$\mu(\omega) \in M \quad \text{for all } \omega \in A_0;$$

see Theorem 12.

Since $\mathcal{B}$ is c.g., $P(\mu \neq \nu) = 0$ for any other r.c.d. $\nu$, and this is basic for (3) and (4). For some time, we guessed that $\mathcal{B}$ c.g. is enough for (3) and (4). Instead, as we now prove, some extra conditions are needed. Let

$$\mathcal{N} = \{B \in \mathcal{B} : P(B) = 0\}.$$

EXAMPLE 1 [*A failure of condition* (3)]. Let $\Omega = \mathbb{R}$, $\mathcal{B}$ the Borel $\sigma$-field, $Q$ a probability measure on $\mathcal{B}$ vanishing on singletons, and $P = \frac{1}{2}(Q + \delta_0)$. If $\mathcal{A} = \sigma(\mathcal{N})$, then $\mu = P$ a.s. and $H(\omega) = \{\omega\}$ for all $\omega$, so that $H(0) = \{0\} \notin \mathcal{A}$ and $\mu(0)\{0\} = P\{0\} = \frac{1}{2}$.

Incidentally, Example 1 exhibits also a couple of (perhaps unexpected) facts. Unless $H(\omega) \in \mathcal{A}$ for all $\omega \in A_0$, (4) does not imply (3). Further, it may be that $\mu(\omega)(H(\omega)) < 1$, for a single point $\omega$ ($\omega = 0$ in Example 1), even though $H(\omega) \in \mathcal{B}$ and $\mu(\omega)(A) = I_A(\omega)$ for all $A \in \mathcal{A}$.

EXAMPLE 2 [*A failure of condition* (4)]. Let $\Omega = \mathbb{R}^2$, $\mathcal{B}$ the Borel $\sigma$-field, and $P = Q \times Q$ where $Q$ is the $N(0,1)$ law on the real Borel sets. Denoting $\mathcal{G}$ the $\sigma$-field on $\Omega$ generated by $(x,y) \mapsto x$, a (natural) r.c.d. for $P$ given $\mathcal{G}$ is $\mu((x,y)) = \delta_x \times Q$. With such a $\mu$, condition (4) fails if $\mathcal{A}$ is taken to be $\mathcal{A} = \sigma(\mathcal{G} \cup \mathcal{N})$. In fact, $\mu$ is also a r.c.d. for $P$ given $\mathcal{A}$, and for all $(x,y)$ one has $\{x\} \times [0, \infty) \in \mathcal{A}$ while

$$\mu((x,y))(\{x\} \times [0, \infty)) = \tfrac{1}{2}.$$

Note also that (3) holds in this example, since $H(\omega) = \{\omega\}$ and $\mu(\omega)\{\omega\} = 0$ for all $\omega \in \Omega$.

This paper provides necessary and sufficient conditions for (3) and (4). Special attention is devoted to the particular case where $\mathcal{A}$ is a *tail sub-$\sigma$-field*, that is, $\mathcal{A}$ is the intersection of a nonincreasing sequence of *countably generated* sub-$\sigma$-fields. The main results are Theorems 3, 4, 8 and 15. Theorem 15 states that (4) is always true whenever $\mathcal{A}$ is a tail sub-$\sigma$-field. Theorems 3, 4 and 8 deal with condition (3). One consequence of Theorem 4 is that, when $(\Omega, \mathcal{B})$ is a standard space, (3) holds for various classically interesting sub-$\sigma$-fields $\mathcal{A}$, including tail, symmetric, invariant, as well as some sub-$\sigma$-fields connected with continuous time processes.

**2. When regular conditional distributions are 0–1 on the (appropriate) atoms of the conditioning $\sigma$-field.** This section deals with condition (3). It is split into three subsections.

2.1. *Basic results.* For any map $\nu : \Omega \to \mathbb{P}$, we write $\sigma(\nu)$ for the $\sigma$-field generated by $\nu(B)$ for all $B \in \mathcal{B}$, where $\nu(B)$ stands for the real function $\omega \mapsto \nu(\omega)(B)$. Since $\mathcal{B}$ is c.g., $\sigma(\nu)$ is c.g., too. In particular, $\sigma(\mu)$ is c.g. Let

$$\mathcal{A}_P = \{A \subset \Omega : \exists A_1, A_2 \in \mathcal{A} \text{ with } A_1 \subset A \subset A_2 \text{ and } P(A_2 - A_1) = 0\}$$

be the completion of $\mathcal{A}$ with respect to $P|\mathcal{A}$. The only probability measure on $\mathcal{A}_P$ agreeing with $P$ on $\mathcal{A}$ is still denoted by $P$. Further, in case $H(\omega) \in \mathcal{B}$ for all $\omega$, we let

$$f_B(\omega) = \mu(\omega)(B \cap H(\omega)) \qquad \text{for } \omega \in \Omega \text{ and } B \in \mathcal{B},$$

$$f = f_\Omega, \qquad S = \{f > 0\}.$$

We are in a position to state our first characterization of (3).

THEOREM 3. *Suppose* (1) *holds and* $H(\omega) \in \mathcal{B}$ *for all* $\omega$. *For each* $U \in \mathcal{A}$ *such that the trace $\sigma$-field* $\mathcal{A} \cap U$ *is c.g., there is* $U_0 \in \mathcal{A}$ *with* $U_0 \subset U$, $P(U - U_0) = 0$ *and* $f(\omega) = 1$ *for all* $\omega \in U_0$. *Moreover, if* $S \in \mathcal{A}_P$, *then condition* (3) *is equivalent to each of the following conditions* (a)–(b):

(a) $f_B$ *is* $\mathcal{A}_P$-*measurable for all* $B \in \mathcal{B}$;
(b) $\mathcal{A} \cap U$ *is c.g. for some* $U \in \mathcal{A}$ *with* $U \subset S$ *and* $P(S - U) = 0$.

PROOF. Suppose $\mathcal{A} \cap U$ is c.g. for some $U \in \mathcal{A}$ and define

$$\mathcal{A}_0 = \{(A \cap U) \cup F : A \in \mathcal{A}, F = \varnothing \text{ or } F = U^c\}.$$

Let $H_0(\omega)$ be the $\mathcal{A}_0$-atom including $\omega$. Then, $\mathcal{A}_0$ is c.g. and $H_0(\omega) = H(\omega)$ for $\omega \in U$. A r.c.d. for $P$ given $\mathcal{A}_0$ is $\mu_0(\omega) = I_U(\omega)\mu(\omega) + I_{U^c}(\omega)\alpha$, where $\alpha(\cdot) = P(\cdot \mid U^c)$ if $P(U) < 1$ and $\alpha$ is any fixed element of $\mathbb{P}$ if $P(U) = 1$. Since $\mathcal{A}_0$ is c.g., there is $K \in \mathcal{A}_0$ with $P(K) = 1$ and $\mu_0(\omega)(H_0(\omega)) = 1$ for all $\omega \in K$. Since $f(\omega) = \mu(\omega)(H(\omega)) = \mu_0(\omega)(H_0(\omega)) = 1$ for each $\omega \in K \cap U$, it suffices to let $U_0 = K \cap U$.

Next, suppose $S \in \mathcal{A}_P$ and take $C, D \in \mathcal{A}$ such that $C \subset S$, $D \subset S^c$ and $P(C \cup D) = 1$.

"(3) $\Rightarrow$ (a)." Let $A = A_0 \cap C$, where $A_0 \in \mathcal{A}$ is such that $P(A_0) = 1$ and $f(\omega) \in \{0, 1\}$ for all $\omega \in A_0$. Fix $B \in \mathcal{B}$. Since $f_B \leq f = 0$ on $D$ and $f_B = \mu(B)$ on $A$, one obtains $I_{A \cup D} f_B = I_A f_B = I_A \mu(B)$. Thus, $f_B$ is $\mathcal{A}_P$-measurable.

"(a) $\Rightarrow$ (b)." Given any $\alpha \in \mathbb{P}$, define the map $\nu : \Omega \to \mathbb{P}$ by

$$\nu(\omega)(B) = I_S(\omega)\frac{f_B(\omega)}{f(\omega)} + I_{S^c}(\omega)\alpha(B), \qquad \omega \in \Omega, B \in \mathcal{B}.$$

Then, (a) implies that $\sigma(\nu) \subset \mathcal{A}_P$. Fix a countable field $\mathcal{B}_0$ generating $\mathcal{B}$, and, for each $B \in \mathcal{B}_0$, take a set $A_B \in \mathcal{A}$ such that $P(A_B) = 1$ and $I_{A_B}\nu(B)$ is $\mathcal{A}$-measurable. Define $U = (\bigcap_{B \in \mathcal{B}_0} A_B) \cap C$ and note that $U \in \mathcal{A}$, $U \subset S$ and $P(S - U) = 0$. Since $I_U \nu(B)$ is $\mathcal{A}$-measurable for each $B \in \mathcal{B}_0$, it follows that $\sigma(\nu) \cap U \subset \mathcal{A} \cap U$. Since $A \cap U = \{\nu(A) = 1\} \cap U$ for all $A \in \mathcal{A}$, then $\mathcal{A} \cap U \subset \sigma(\nu) \cap U$. Hence, $\mathcal{A} \cap U = \sigma(\nu) \cap U$ is c.g.

"(b) $\Rightarrow$ (3)." By the first assertion of the theorem, since $U \in \mathcal{A}$ and $\mathcal{A} \cap U$ is c.g., there is $U_0 \in \mathcal{A}$ with $U_0 \subset U$, $P(U - U_0) = 0$ and $f = 1$ on $U_0$. Define $A_0 = U_0 \cup D$ and note that $A_0 \in \mathcal{A}$ and $f \in \{0, 1\}$ on $A_0$. Since $U \subset S$ and $P(S - U) = 0$, one also obtains $P(A_0) = P(U_0) + P(D) = P(S) + P(S^c) = 1$. □

A basic condition for existence of disintegrations is that

$$G = \{(x, y) \in \Omega \times \Omega : H(x) = H(y)\}$$

belongs to $\mathcal{B} \otimes \mathcal{B}$; see [1]. Such a condition also plays a role in our main characterization of (3). Let

$$\mathcal{A}^* = \{B \in \mathcal{B} : B \text{ is a union of } \mathcal{A}\text{-atoms}\}.$$

THEOREM 4. *If (1) holds and $G \in \mathcal{B} \otimes \mathcal{B}$, then $H(\omega) \in \mathcal{B}$ for all $\omega \in \Omega$ and $f_B$ is $\mathcal{A}^*$-measurable for all $B \in \mathcal{B}$. If in addition $S \in \mathcal{A}_P$, then each of conditions (3), (a) and (b) is equivalent to*

(c) $\mathcal{A} \cap U = \mathcal{A}^* \cap U$ *for some $U \in \mathcal{A}$ with $U \subset S$ and $P(S - U) = 0$.*

PROOF. For $C \subset \Omega \times \Omega$, let $C_\omega = \{u \in \Omega : (\omega, u) \in C\}$ be the $\omega$-section of $C$. Suppose (1) holds and $G \in \mathcal{B} \otimes \mathcal{B}$. Then, $H(\omega) = G_\omega \in \mathcal{B}$ for all $\omega$. By a monotone class argument, the map $\omega \mapsto \mu(\omega)(C_\omega)$ is $\mathcal{B}$-measurable whenever $C \in \mathcal{B} \otimes \mathcal{B}$. Letting $C = G \cap (\Omega \times B)$, where $B \in \mathcal{B}$, implies that $f_B$ is $\mathcal{B}$-measurable. Since $f_B$ is constant on each $\mathcal{A}$-atom, it is in fact $\mathcal{A}^*$-measurable. (Note that $\mathcal{A}^*$-measurability of $f_B$ does not require $\mathcal{B}$ c.g.) Next, suppose also that $S \in \mathcal{A}_P$. By Theorem 3, conditions (3), (a) and (b) are equivalent. Suppose (c) holds, and fix $B \in \mathcal{B}$ and a Borel set $I \subset \mathbb{R}$. Since $\{f_B \in I\} \in \mathcal{A}^*$, condition (c) yields $\{f_B \in I\} \cap U \in \mathcal{A}$, and $\{f_B \in I\} \cap (S - U) \in \mathcal{A}_P$ due to $(S - U) \in \mathcal{A}_P$ and $P(S - U) = 0$. Thus (assuming $0 \in I$ to fix ideas),

$$\{f_B \in I\} = S^c \cup (\{f_B \in I\} \cap (S - U)) \cup (\{f_B \in I\} \cap U) \in \mathcal{A}_P,$$

so that (a) holds. Conversely, suppose (a) holds. For each $B \in \mathcal{B}$, there is $A_B \in \mathcal{A}$ such that $P(A_B) = 1$ and $I_{A_B} f_B$ is $\mathcal{A}$-measurable. Letting $A = \bigcap_{B \in \mathcal{B}_0} A_B$, where $\mathcal{B}_0$ is a countable field generating $\mathcal{B}$, it follows that $A \in \mathcal{A}$, $P(A) = 1$ and $I_A f_B$ is $\mathcal{A}$-measurable for all $B \in \mathcal{B}$. Since $S \in \mathcal{A}_P$ and $P(A) = 1$, there is $U \in \mathcal{A}$ with $U \subset A \cap S$ and $P(S - U) = 0$. Given $B \in \mathcal{A}^*$, on noting that $f_B = I_B f$, one obtains

$$B \cap U = \{I_B f > 0\} \cap U = \{f_B > 0\} \cap U = (\{f_B > 0\} \cap A) \cap U \in \mathcal{A}.$$

Hence $\mathcal{A} \cap U = \mathcal{A}^* \cap U$, that is, condition (c) holds. □

We now state a couple of corollaries to Theorem 4. The first covers in particular the case where the $\mathcal{A}$-atoms are the singletons, while the second (and more important) applies to various real situations.

COROLLARY 5. *Suppose* (1) *holds,* $S \in \mathcal{A}_P$ *and* $\mathcal{A}$, $\mathcal{B}$ *have the same atoms. Then,* (3) *holds if and only if* $\mathcal{A} \cap U = \mathcal{B} \cap U$ *for some* $U \in \mathcal{A}$ *with* $U \subset S$ *and* $P(S - U) = 0$.

PROOF. Since $\mathcal{A}$, $\mathcal{B}$ have the same atoms and $\mathcal{B}$ is c.g.,

$$G = \{(x, y) : x \text{ and } y \text{ are in the same } \mathcal{B}\text{-atom}\} \in \mathcal{B} \otimes \mathcal{B}.$$

Therefore, it suffices applying Theorem 4 and noting that $\mathcal{A}^* = \mathcal{B}$. □

COROLLARY 6. *If* (1) *holds,* $G \in \mathcal{B} \otimes \mathcal{B}$ *and* $\mathcal{A} \cap C = \mathcal{A}^* \cap C$, *for some* $C \in \mathcal{A}$ *with* $P(C) = 1$, *then condition* (3) *holds.*

PROOF. Since $C \in \mathcal{A}$ and $S \cap C \in \mathcal{A}^* \cap C = \mathcal{A} \cap C$, then $S \cap C \in \mathcal{A}$. Since $P(C) = 1$, it follows that $S \in \mathcal{A}_P$ and (c) holds with $U = S \cap C$. Thus, (3) follows from Theorem 4. □

As shown in [4], if $(\Omega, \mathcal{B})$ is a standard space ($\Omega$ Borel subset of a Polish space and $\mathcal{B}$ the Borel $\sigma$-field on $\Omega$), then $G \in \mathcal{B} \otimes \mathcal{B}$ and $\mathcal{A}^* = \mathcal{A}$ for various classically interesting sub-$\sigma$-fields $\mathcal{A}$, including tail, symmetric, invariant, as well as some sub-$\sigma$-fields connected with continuous time processes. In view of Corollary 6, *condition* (3) *holds in case* $(\Omega, \mathcal{B})$ *is a standard space and* $\mathcal{A}$ *is any one of the above-mentioned sub-$\sigma$-fields*.

2.2. *Tail sub-$\sigma$-fields.* When condition (3) holds, the next step is determining those $\omega$'s satisfying $f(\omega) = 1$. Suppose the assumptions of Corollary 6 are in force [so that (3) holds and $f$ is $\mathcal{A}_P$-measurable] and define

$$\mathcal{U} = \{U \in \mathcal{A} : \mathcal{A} \cap U \text{ is c.g.}\} \cup \{\varnothing\}.$$

Since $\mathcal{U}$ is closed under countable unions, some $A \in \mathcal{U}$ meets $P(A) = \sup\{P(U) : U \in \mathcal{U}\}$. By the first assertion in Theorem 3, $P(A - \{f = 1\}) = 0$. Taking $U$ as in condition (b) and noting that $U \in \mathcal{U}$, one also obtains $P(\{f = 1\} - A) = P(U - A) = 0$. Therefore, $A$ is the set we are looking for, in the sense that

$$P(\{f = 1\} \Delta A) = 0.$$

Incidentally, the above remarks provide also a criterion for deciding whether $\mu$ is *maximally improper* according to [11]. Under the assumptions of Corollary 6, in fact, $\mu$ is maximally improper precisely when $P(S) = 0$. Hence,

$$\mu \text{ is maximally improper} \quad \Leftrightarrow \quad P(U) = 0 \text{ for all } U \in \mathcal{U}.$$

Some handy description of the members of $\mathcal{U}$, thus, would be useful. Unfortunately, such a description is generally hard to be found. We now discuss a particular case.

Let $\mathcal{A}$ be a *tail sub-$\sigma$-field*, that is, $\mathcal{A} = \bigcap_{n \geq 1} \mathcal{A}_n$ where $\mathcal{A}_n$ is a *countably generated* $\sigma$-field and $\mathcal{B} \supset \mathcal{A}_n \supset \mathcal{A}_{n+1}$ for all $n \geq 1$. As already noted, the assumptions of Corollary 6 hold for such an $\mathcal{A}$ if $(\Omega, \mathcal{B})$ is a standard space. More generally, it is enough that:

LEMMA 7. *If $\mathcal{A}$ is a tail sub-$\sigma$-field, (1) holds and*

*for each $n$, there is a r.c.d. $\mu_n$ for $P$ given $\mathcal{A}_n$,*

*then $G \in \mathcal{B} \otimes \mathcal{B}$ and $\mathcal{A} \cap C = \mathcal{A}^* \cap C$ for some $C \in \mathcal{A}$ with $P(C) = 1$.*

PROOF. Since $G_n := \{(x,y) \in \Omega \times \Omega : x \text{ and } y \text{ are in the same } \mathcal{A}_n\text{-atom}\} \in \mathcal{A}_n \otimes \mathcal{A}_n$, Proposition 1 of [4] implies $G = \bigcup_n G_n \in \mathcal{B} \otimes \mathcal{B}$. For each $n$, since $\mathcal{A}_n$ is c.g., there is $C_n \in \mathcal{A}_n$ such that $P(C_n) = 1$ and $\mu_n(\omega)(A) = I_A(\omega)$ whenever $A \in \mathcal{A}_n$ and $\omega \in C_n$. Define $C = \bigcup_{n \geq 1} \bigcap_{j \geq n} C_j$ and note that $C \in \mathcal{A}$ and $P(C) = 1$. Fix $B \in \mathcal{A}^*$. Since $B$ is a union of $\mathcal{A}_n$-atoms whatever $n$ is,

$$\lim_n \mu_n(\omega)(B) = I_B(\omega) \quad \text{for all } \omega \in C.$$

Thus, $B \cap C = \{\lim_n \mu_n(B) = 1\} \cap C \in \mathcal{A} \cap C$. □

Each $\mathcal{A}_n$, being c.g., can be written as $\mathcal{A}_n = \sigma(X_n)$ for some $X_n : \Omega \to \mathbb{R}$. Since $\mathcal{A}_n \supset \mathcal{A}_j$ for $j \geq n$, it follows that $\mathcal{A}_n = \sigma(X_n, X_{n+1}, \ldots)$. Thus, $\mathcal{A}$ admits the usual representation

$$\mathcal{A} = \bigcap_n \sigma(X_n, X_{n+1}, \ldots)$$

for some sequence $(X_n)$ of real random variables. In particular,

$$H(\omega) = \{\exists n \geq 1 \text{ such that } X_j = X_j(\omega) \text{ for all } j \geq n\} \in \mathcal{A}$$

so that $\mathcal{A}$ includes its atoms. Note also that a c.g. sub-$\sigma$-field is tail while the converse need not be true. In fact, for a $\sigma$-field $\mathcal{F}$ to be not c.g., it is enough that $\mathcal{F}$ supports a 0–1 valued probability measure $Q$ such that $Q(F) = 0$ whenever $F \in \mathcal{F}$ and $F$ is an $\mathcal{F}$-atom; see Theorem 1 of [4]. Thus, for instance, $\mathcal{A} = \bigcap_n \sigma(X_n, X_{n+1}, \ldots)$ is not c.g. in case $(X_n)$ is i.i.d. and $X_1$ has a nondegenerate distribution.

To find usable characterizations of $\mathcal{U}$ is not an easy task. Countable unions of $\mathcal{A}$-atoms belong to $\mathcal{U}$, but generally they are not all the elements of $\mathcal{U}$. For instance, if $\Omega = \mathbb{R}^\infty$ and $X_n$ is the $n$th coordinate projection, then

$$U = \{\exists n \geq 1 \text{ such that } X_j = X_n \text{ for all } j \geq n\}$$

is an uncountable union of $\mathcal{A}$-atoms. However, $\mathcal{A} \cap U$ is c.g. since $\mathcal{A} \cap U = \sigma(L) \cap U$ where $L = \limsup_n X_n$.

Another possibility could be selecting a subclass $\mathbb{Q} \subset \mathbb{P}$ and showing that $U \in \mathcal{U}$ if and only if $U \in \mathcal{A}$ and $Q(U) = 0$ for each $Q \in \mathbb{Q}$. We do not know whether some (nontrivial) characterization of this type is available. Here, we just note that

$$\mathbb{Q}_0 = \{Q \in \mathbb{P} : (X_n) \text{ is i.i.d. and } X_1 \text{ has a nondegenerate distribution, under } Q\}$$

does not work (though the "only if" implication is true, in view of Theorem 1 of [4]). As an example, $U := \{X_n \to 0\} \notin \mathcal{U}$ even though $U \in \mathcal{A}$ and $Q(U) = 0$ for all $Q \in \mathbb{Q}_0$. To see that $U \notin \mathcal{U}$, let $X_n$ be the $n$th coordinate projection on $\Omega = \mathbb{R}^\infty$, and let $P_U$ be a probability measure on the Borel sets of $\Omega$ which makes $(X_n)$ independent and each $X_n$ uniformly distributed on $(0, \frac{1}{n})$. Then $P_U(U) = 1$ and, when restricted to $\mathcal{A} \cap U$, $P_U$ is a 0–1 probability measure such that $P_U(H(\omega)) = 0$ for each $\omega \in U$. Hence, Theorem 1 of [4] implies that $\mathcal{A} \cap U$ is not c.g.

A last note is that $P(S)$ can assume any value between 0 and 1. For instance, take $U \in \mathcal{U}$ and $P_1, P_2 \in \mathbb{P}$ such that: (i) $P_1(U) = P_2(U^c) = 1$; (ii) $P_2$ is 0–1 on $\mathcal{A}$ with $P_2(H(\omega)) = 0$ for all $\omega$. Define $P = uP_1 + (1-u)P_2$ where $u \in (0,1)$. A r.c.d. for $P$ given $\mathcal{A}$ is $\mu(\omega) = I_U(\omega)\mu_1(\omega) + I_{U^c}(\omega)P_2$, where $\mu_1$ denotes a r.c.d. for $P_1$ given $\mathcal{A}$. Since $U \in \mathcal{U}$, Theorem 3 implies $\mu_1(\omega)(H(\omega)) = 1$ for $P_1$-almost all $\omega \in U$. Thus, $P(S) = P(U) = u$.

2.3. *Miscellaneous results.* A weaker version of (3) lies in asking $\mu(\cdot)(H(\cdot))$ to be 0–1 over a set of $\mathcal{A}^*$, but not necessarily of $\mathcal{A}$, that is

(3\*) There is $B_0 \in \mathcal{A}^*$ with $P(B_0) = 1$ and $\mu(\omega)(H(\omega)) \in \{0, 1\}$ for all $\omega \in B_0$. Suitably adapted, the proofs of Theorems 3 and 4 yield a characterization of (3\*) as well. Recall $\mathcal{N} = \{B \in \mathcal{B} : P(B) = 0\}$ and note that

$$\sigma(\mathcal{A} \cup \mathcal{N}) = \{B \in \mathcal{B} : \mu(B) = I_B \text{ a.s.}\}.$$

THEOREM 8. *Suppose* (1) *holds and* $G \in \mathcal{B} \otimes \mathcal{B}$. *Then, condition* (3\*) *implies* $S \in \sigma(\mathcal{A} \cup \mathcal{N})$. *Moreover, if* $S \in \sigma(\mathcal{A} \cup \mathcal{N})$, *then*

$$(3*) \quad \Leftrightarrow \quad (b*) \quad \Leftrightarrow \quad (c*)$$

*where:*

(b\*) $\mathcal{A} \cap V$ *is c.g. for some* $V \in \mathcal{A}^*$ *with* $V \subset S$ *and* $P(S - V) = 0$;
(c\*) $\mathcal{A} \cap V = \mathcal{A}^* \cap V$ *for some* $V \in \mathcal{A}^*$ *with* $V \subset S$ *and* $P(S - V) = 0$.

PROOF. If (3*) holds, then $\mu(S) = 1$ on $B_0 \cap S$, and since $P(B_0) = 1$ one obtains

$$E(\mu(S)I_{S^c}) = P(S) - E(\mu(S)I_{B_0}I_S) = P(S) - E(I_{B_0}I_S) = 0.$$

Thus, $\mu(S) = I_S$ a.s., that is, $S \in \sigma(\mathcal{A} \cup \mathcal{N})$. Next, suppose that $S \in \sigma(\mathcal{A} \cup \mathcal{N})$.

"(3*) $\Rightarrow$ (c*)." Define $V = B_0 \cap S$ and note that $B \cap V = \{\mu(B) = 1\} \cap V \in \mathcal{A} \cap V$ for all $B \in \mathcal{A}^*$.

"(c*) $\Rightarrow$ (b*)." Fix $\alpha \in \mathbb{P}$ and define $\nu(\omega)(B) = I_V(\omega)\frac{f_B(\omega)}{f(\omega)} + I_{V^c}(\omega)\alpha(B)$ for all $\omega \in \Omega$ and $B \in \mathcal{B}$. Then, $\sigma(\nu) \subset \mathcal{A}^*$. Further, $\nu(\omega)(H(\omega)) = 1$ for all $\omega \in V$, so that $B \cap V = \{\nu(B) = 1\} \cap V$ for all $B \in \mathcal{A}^*$. Hence, (c*) implies that $\mathcal{A} \cap V = \mathcal{A}^* \cap V = \sigma(\nu) \cap V$ is c.g.

"(b*) $\Rightarrow$ (3*)." Let $\mathcal{A}_0 = \{(A \cap V) \cup F : A \in \mathcal{A}, F = \varnothing \text{ or } F = V^c\}$. Since $\mu(V) = \mu(S) = I_S = I_V$ a.s., for all $A \in \mathcal{A}$ and $B \in \mathcal{B}$ one obtains

$$E(I_A I_V \mu(B)) = E(I_A \mu(B \cap V)) = P((A \cap V) \cap B).$$

So, $\mu_0(\omega) = I_V(\omega)\mu(\omega) + I_{V^c}(\omega)\alpha$ is a r.c.d. for $P$ given $\mathcal{A}_0$, where $\alpha(\cdot) = P(\cdot \mid V^c)$ if $P(V) < 1$ and $\alpha$ is any fixed element of $\mathbb{P}$ if $P(V) = 1$. Since $\mathcal{A}_0$ is c.g., there is $K \in \mathcal{A}_0$ with $P(K) = 1$ and $\mu_0(\omega)(H_0(\omega)) = 1$ for all $\omega \in K$, where $H_0(\omega)$ denotes the $\mathcal{A}_0$-atom including $\omega$. Hence, it suffices to let $B_0 = (K \cap V) \cup S^c$ and noting that $H_0(\omega) = H(\omega)$ and $\mu_0(\omega) = \mu(\omega)$ for all $\omega \in V$. □

One consequence of Theorem 8 is that, if (1) holds and $G \in \mathcal{B} \otimes \mathcal{B}$, then condition (3*) is equivalent to $\mu(S) = I_S$ a.s. and $P(0 < f \leq \frac{1}{2}) = 0$. In fact,

$$A \cap \{f > \tfrac{1}{2}\} = \{\mu(A) > \tfrac{1}{2}\} \cap \{f > \tfrac{1}{2}\} \qquad \text{for all } A \in \mathcal{A},$$

so that $\mathcal{A} \cap \{f > \tfrac{1}{2}\} = \sigma(\mu) \cap \{f > \tfrac{1}{2}\}$ is c.g. Hence, if $P(0 < f \leq \tfrac{1}{2}) = 0$, condition (b*) holds with $V = \{f > \tfrac{1}{2}\}$.

Finally, we give one more condition for (3). Though seemingly simple, it is hard to be tested in real problems.

PROPOSITION 9. *If (1) holds and $H(\omega) \in \mathcal{B}$ for all $\omega$, a sufficient condition for (3) is*

(5) $$\mu(x)(H(y)) = 0 \qquad \text{whenever } H(x) \neq H(y).$$

PROOF. As stated in the forthcoming Lemma 10, since $\sigma(\mu)$ is c.g. and $\mu$ is also a r.c.d. for $P$ given $\sigma(\mu)$, there is a set $T \in \sigma(\mu)$ such that $P(T) = 1$ and $\mu(\omega)(\mu = \mu(\omega)) = 1$ for each $\omega \in T$. Let $A_0 = T$ and fix $\omega \in S$. Then, $\mu(\omega) = \mu(x)$ if $x \in H(\omega)$ [since $\sigma(\mu) \subset \mathcal{A}$] and $\mu(\omega) \neq \mu(x)$ if $x \notin H(\omega)$ since in the latter case (5) yields

$$\mu(x)(H(\omega)) = 0 < f(\omega) = \mu(\omega)(H(\omega)).$$

Thus, $H(\omega) = \{\mu = \mu(\omega)\}$. If $\omega \in T \cap S = A_0 \cap S$, this implies

$$\mu(\omega)(H(\omega)) = \mu(\omega)(\mu = \mu(\omega)) = 1. \qquad \square$$

**3. When regular conditional distributions are 0–1 on the conditioning $\sigma$-field.** In this section, condition (4) is shown to be true whenever $\mathcal{A}$ is a tail sub-$\sigma$-field. Moreover, two characterizations of (4) and a result in the negative [i.e., a condition for (4) to be false] are given.

We begin by recalling a few simple facts about $\sigma(\mu)$.

LEMMA 10. *If* (1) *holds, then* $\sigma(\mu)$ *is c.g.,* $\mu$ *is a r.c.d. for* $P$ *given* $\sigma(\mu)$, *and there is a set* $T \in \sigma(\mu)$ *with* $P(T) = 1$ *and*

$$\mu(\omega)(\mu = \mu(\omega)) = 1 \qquad \text{for all } \omega \in T.$$

*Moreover,*

$$\mathcal{A} = \sigma(\sigma(\mu) \cup (\mathcal{A} \cap \mathcal{N})).$$

PROOF. Since $\sigma(\mu) \subset \mathcal{A}$, $\mu$ is a r.c.d. given $\sigma(\mu)$. Since $\mathcal{B}$ is c.g., $\sigma(\mu)$ is c.g. with atoms of the form $\{\mu = \mu(\omega)\}$. Hence, there is $T \in \sigma(\mu)$ with $P(T) = 1$ and $\mu(\omega)(\mu = \mu(\omega)) = 1$ for all $\omega \in T$. Finally, since

$$A = (\{\mu(A) = 1\} \cap \{\mu(A) = I_A\}) \cup (A \cap \{\mu(A) \neq I_A\})$$

for all $A \in \mathcal{A}$, it follows that $\mathcal{A} \subset \sigma(\sigma(\mu) \cup (\mathcal{A} \cap \mathcal{N})) \subset \mathcal{A}$. $\square$

By Lemma 10, $\mu(\omega)$ is 0–1 on $\sigma(\mu)$ for each $\omega \in T$. Since $\mathcal{A} = \sigma(\sigma(\mu) \cup (\mathcal{A} \cap \mathcal{N}))$, condition (4) can be written as

$$\mu(\omega)(A) \in \{0, 1\} \qquad \text{for all } \omega \in A_0 \text{ and } A \in \mathcal{A} \text{ with } P(A) = 0.$$

In particular, (4) *holds whenever* $P$ *is atomic on* $\mathcal{A}$, in the sense that there is a countable partition $\{A_1, A_2, \ldots\}$ of $\Omega$ satisfying $A_j \in \mathcal{A}$ and $P(A \cap A_j) \in \{0, P(A_j)\}$ for all $j \geq 1$ and $A \in \mathcal{A}$. In this case, in fact, $\mu(\omega) \ll P$ for each $\omega$ in some set $C \in \mathcal{A}$ with $P(C) = 1$.

Slightly developing the idea underlying Example 2, we next give a sufficient condition for (4) to be false.

PROPOSITION 11. *Suppose* (1) *holds and* $P(\mu = \mu(\omega)) = 0$ *for all* $\omega$. *Then,*

$$F = \{\omega : \mu(\omega) \text{ is not 0–1 on } \mathcal{B}\} \quad \text{and} \quad F_0 = \{\omega : \mu(\omega) \text{ is nonatomic on } \mathcal{B}\}$$

*belong to* $\sigma(\mu)$. *Moreover, if* $\mathcal{N} \subset \mathcal{A}$, *then*

$\mu(\omega)$ *is not 0–1 on* $\mathcal{A}$ *for each* $\omega \in F \cap T$, *and*

$\mu(\omega)$ *is nonatomic on* $\mathcal{A}$ *for each* $\omega \in F_0 \cap T$

*with* $T$ *as in Lemma* 10. *In particular, condition* (4) *fails if* $P(F) > 0$.

PROOF. Since $\mathcal{B}$ is c.g., it is clear that $F \in \sigma(\mu)$, while $F_0 \in \sigma(\mu)$ is from [6] (see Corollary 2.13, page 1214). Suppose now that $\mathcal{N} \subset \mathcal{A}$. Let $\omega \in F \cap T$. Since $\omega \in F$, there is $B_\omega \in \mathcal{B}$ with $\mu(\omega)(B_\omega) \in (0,1)$. Define $A_\omega = B_\omega \cap \{\mu = \mu(\omega)\}$. Since $\mathcal{N} \subset \mathcal{A}$ and $P(A_\omega) \leq P(\mu = \mu(\omega)) = 0$, then $A_\omega \in \mathcal{A}$. Since $\omega \in T$,

$$\mu(\omega)(A_\omega) = \mu(\omega)(B_\omega) \in (0,1)$$

so that $\mu(\omega)$ is not 0–1 on $\mathcal{A}$. Finally, fix $\omega \in F_0 \cap T$ and $\varepsilon > 0$. Since $\omega \in F_0$, there is a finite partition $\{B_{1,\omega}, \ldots, B_{n,\omega}\}$ of $\Omega$ such that $B_{i,\omega} \in \mathcal{B}$ and $\mu(\omega)(B_{i,\omega}) < \varepsilon$ for all $i$. As above, letting $A_{i,\omega} = B_{i,\omega} \cap \{\mu = \mu(\omega)\}$, one obtains $A_{i,\omega} \in \mathcal{A}$ and $\mu(\omega)(A_{i,\omega}) = \mu(\omega)(B_{i,\omega}) < \varepsilon$. Hence, $\mu(\omega)$ is nonatomic on $\mathcal{A}$ since $\mu(\omega)(\mu \neq \mu(\omega)) = 0$. □

Even if $\mathcal{N}$ is not contained in $\mathcal{A}$, Proposition 11 applies at least to $\mathcal{A}' = \sigma(\mathcal{A} \cup \mathcal{N})$. Under mild conditions, $\mu$ is even nonatomic on $\mathcal{A}'$ with probability $P(F_0)$. Thus, a lot of r.c.d.'s give rise to a failure of (4) on some sub-$\sigma$-field $\mathcal{A}'$. Since we are conditioning to $\mathcal{A}$ (and not to $\mathcal{A}'$), this fact is not essential. On the other hand, it suggests that (4) is a rather delicate condition.

If $P$ is invariant under a countable collection of measurable transformations and $\mathcal{A}$ is the corresponding invariant sub-$\sigma$-field, then (4) holds; see [9]. This well-known fact is generalized by our first characterization of (4).

THEOREM 12. *Suppose (1) holds and let $M = \{Q \in \mathbb{P} : \mu$ is a r.c.d. for $Q$ given $\mathcal{A}\}$. Then, $Q$ is an extreme point of $M$ if and only if $Q \in M$ and $Q$ is 0–1 on $\mathcal{A}$, and in that case $Q = \mu(\omega)$ for some $\omega \in \Omega$. Moreover, for each $\omega \in T$ (with $T$ as in Lemma 10), the following statements are equivalent:*

(i) $\mu(\omega)(A) \in \{0,1\}$ *for all $A \in \mathcal{A}$;*
(ii) $\mu(\omega)$ *is an extreme point of $M$;*
(iii) $\mu(\omega) \in M$.

*In particular, condition (4) holds if and only if, for some $A_0 \in \mathcal{A}$ with $P(A_0) = 1$,*

$$\mu(\omega) \in M \quad \text{for all } \omega \in A_0.$$

PROOF. Fix $Q \in M$. If $Q(A) \in (0,1)$ for some $A \in \mathcal{A}$, then

$$Q(\cdot) = Q(A)Q(\cdot \mid A) + (1 - Q(A))Q(\cdot \mid A^c),$$

and $Q$ is not extreme since $Q(\cdot \mid A)$ and $Q(\cdot \mid A^c)$ are distinct elements of $M$. Suppose now that $Q = uQ_1 + (1-u)Q_2$, where $u \in (0,1)$ and $Q_1 \neq Q_2$ are in $M$. Since two elements of $M$ coincide if and only if they coincide on $\mathcal{A}$, there is $A \in \mathcal{A}$ with $Q_1(A) \neq Q_2(A)$, and this implies $Q(A) \in (0,1)$. Hence, $Q \in M$ is extreme if and only if it is 0–1 on $\mathcal{A}$. In particular, if $Q$ is extreme then it is 0–1 on the c.g. $\sigma$-field $\sigma(\mu)$, so that $Q(\mu = \mu(\omega)) = 1$ for some $\omega \in \Omega$; see Theorem 1 of [4]. Thus,

$$Q(B) = \int \mu(x)(B) \, Q(dx) = \mu(\omega)(B) \qquad \text{for all } B \in \mathcal{B}.$$

This concludes the proof of the first part. As to the second one, fix $\omega \in T$, and let $A$ and $B$ denote arbitrary elements of $\mathcal{A}$ and $\mathcal{B}$, respectively. Since $\omega \in T$,

$$\int_A \mu(x)(B)\mu(\omega)(dx) = \int_{A \cap \{\mu = \mu(\omega)\}} \mu(x)(B)\mu(\omega)(dx) = \mu(\omega)(A)\mu(\omega)(B).$$

"(i) $\Rightarrow$ (ii)." By what already proved, it is enough showing that $\mu(\omega) \in M$, and this depends on $\mu(\omega)(A \cap B) = \mu(\omega)(A)\mu(\omega)(B) = \int_A \mu(x)(B) \, \mu(\omega)(dx)$.

"(ii) $\Rightarrow$ (iii)." Obvious.

"(iii) $\Rightarrow$ (i)." Under (iii), $\mu(\omega)(A \cap B) = \int_A \mu(x)(B) \, \mu(\omega)(dx) = \mu(\omega) \times (A)\mu(\omega)(B)$, and letting $B = A$ yields $\mu(\omega)(A) = \mu(\omega)(A)^2$. $\square$

Next characterization of (4) stems from a result of [8], Lemma 2A, page 391.

THEOREM 13 (Fremlin). *Let $X$ be an Hausdorff topological space, $\mathcal{F}$ a $\sigma$-field on $X$ including the open sets, $Q$ a complete Radon probability measure on $\mathcal{F}$, and $\mathcal{C}_0$ a class of pairwise disjoint $Q$-null elements of $\mathcal{F}$. Then,*

$$\bigcup_{C \in \mathcal{C}} C \in \mathcal{F} \text{ for all } \mathcal{C} \subset \mathcal{C}_0 \iff Q\left(\bigcup_{C \in \mathcal{C}_0} C\right) = 0.$$

Say that $P$ is *perfect* in case each $\mathcal{B}$-measurable function $h : \Omega \to \mathbb{R}$ meets $P(h \in I) = 1$ for some real Borel set $I \subset h(\Omega)$. For $P$ to be perfect, it is enough that $\Omega$ is an universally measurable subset of a Polish space and $\mathcal{B}$ the Borel $\sigma$-field on $\Omega$. In the present framework, since $\mathcal{B}$ is c.g., Theorem 13 applies precisely when $P$ is perfect. We are now able to state our second characterization of (4). It is of possible theoretical interest even if of little practical use.

THEOREM 14. *Suppose* (1) *holds and $P$ is perfect, define*

$$\mathcal{A}(\omega) = \{A \in \mathcal{A} : \mu(\omega)(A) \in \{0,1\}\} \qquad \omega \in \Omega,$$

*and let $\Gamma_0$ denote the class of those $\sigma$-fields $\mathcal{G} \subset \mathcal{A}$ with $\mathcal{G} \neq \mathcal{A}$. Then, condition* (4) *holds if and only if*

(6) $$\bigcup_{\mathcal{G} \in \Gamma} \{\omega : \mathcal{A}(\omega) = \mathcal{G}\} \in \mathcal{A}_P \qquad \text{for all } \Gamma \subset \Gamma_0.$$

PROOF. If $\mu(\omega)$ is 0–1 on $\mathcal{A}$ for all $\omega \in A_0$, where $A_0 \in \mathcal{A}$ and $P(A_0) = 1$, then (6) follows from

$$\{\omega : \mathcal{A}(\omega) = \mathcal{G}\} \subset A_0^c \qquad \text{for all } \mathcal{G} \in \Gamma_0.$$

Conversely, suppose (6) holds. Let $X$ be the partition of $\Omega$ in the atoms of $\mathcal{B}$. The elements of $\mathcal{B}$ are unions of elements of $X$, so that $\mathcal{B}$ can be regarded as a $\sigma$-field on $X$. Let $(X, \mathcal{F}, Q)$ be the completion of $(X, \mathcal{B}, P)$. Since $\mathcal{B}$ is c.g., under a suitable distance, $X$ is separable metric and $\mathcal{B}$ the corresponding Borel $\sigma$-field; see [2]. Since $P$ is perfect, $P$ is Radon by a result of Sazonov (Theorem 12 of [10]), so that $Q$ is Radon, too. Next, define $C_\mathcal{G} = \{\omega : \mathcal{A}(\omega) = \mathcal{G}\}$ for $\mathcal{G} \in \Gamma_0$, $U_A = \{\omega : \mu(\omega)(A) \in (0, 1)\}$ for $A \in \mathcal{A}$, and $U = \{\omega : \mathcal{A}(\omega) \neq \mathcal{A}\}$ (all regarded as subsets of $X$). For each $\mathcal{G} \in \Gamma_0$ there is $A \in \mathcal{A}$ with $C_\mathcal{G} \subset U_A$. Since $U_A \in \mathcal{A}$ and $P(U_A) = 0$, then $C_\mathcal{G} \in \mathcal{F}$ and $Q(C_\mathcal{G}) = 0$. Hence, $\mathcal{C}_0 = \{C_\mathcal{G} : \mathcal{G} \in \Gamma_0\}$ is a collection of pairwise disjoint $Q$-null elements of $\mathcal{F}$ satisfying $U = \bigcup_{\mathcal{G} \in \Gamma_0} C_\mathcal{G}$. By (6), Theorem 13 yields $Q(U) = 0$. Finally, since $U \in \mathcal{A}_P$, $Q(U) = 0$ implies $U \subset A$ for some $A \in \mathcal{A}$ with $P(A) = 0$. Thus, to get (4), it suffices to let $A_0 = A^c$. □

Finally, by a martingale argument, we prove that (4) holds when $\mathcal{A}$ is a tail sub-$\sigma$-field. This is true, in addition, even though $\mathcal{B}$ fails to be c.g.

THEOREM 15. *Let $\mathcal{A} = \bigcap_{n \geq 1} \mathcal{A}_n$, where $\mathcal{B} \supset \mathcal{A}_1 \supset \mathcal{A}_2 \supset \cdots$ and $\mathcal{A}_n$ is a c.g. $\sigma$-field for each $n$. Given a r.c.d. $\mu$, for $P$ given $\mathcal{A}$, there is a set $A_0 \in \mathcal{A}$ such that $P(A_0) = 1$ and $\mu(\omega)(A) \in \{0, 1\}$ for all $A \in \mathcal{A}$ and $\omega \in A_0$.*

PROOF. First recall that a probability measure $Q \in \mathbb{P}$ is 0–1 on $\mathcal{A}$ if (and only if) $\sup_{A \in \mathcal{A}_n} |Q(A \cap B) - Q(A)Q(B)| \to 0$, as $n \to \infty$, for all $B \in \mathcal{A}_1$. Also, given any field $\mathcal{F}_n$ such that $\mathcal{A}_n = \sigma(\mathcal{F}_n)$, the "sup" can be taken over $\mathcal{F}_n$, that is,

$$\sup_{A \in \mathcal{A}_n} |Q(A \cap B) - Q(A)Q(B)| = \sup_{A \in \mathcal{F}_n} |Q(A \cap B) - Q(A)Q(B)|.$$

Now, since the $\mathcal{A}_n$ are c.g., there are countable fields $\mathcal{F}_n$ satisfying $\mathcal{A}_n = \sigma(\mathcal{F}_n)$ for all $n$. Let

$$V_n^B(\omega) = \sup_{A \in \mathcal{F}_n} |\mu(\omega)(A \cap B) - \mu(\omega)(A)\mu(\omega)(B)|, \qquad n \geq 1, B \in \mathcal{A}_1, \omega \in \Omega.$$

Since $\mathcal{F}_n$ is countable, $V_n^B$ is an $\mathcal{A}$-measurable random variable for all $n$ and $B$. It is enough proving that

(7) $\qquad V_n^B \to 0 \qquad \text{a.s., as } n \to \infty, \text{ for all } B \in \mathcal{A}_1.$

Suppose in fact (7) holds and define

$$A_0 = \left\{ \omega : \lim_n V_n^B(\omega) = 0 \text{ for each } B \in \mathcal{F}_1 \right\}.$$

Since $\mathcal{F}_1$ is countable, $A_0 \in \mathcal{A}$ and (7) implies $P(A_0) = 1$. Fix $\omega \in A_0$. Since $\mathcal{A}_1 = \sigma(\mathcal{F}_1)$, given $B \in \mathcal{A}_1$ and $\varepsilon > 0$, there is $C \in \mathcal{F}_1$ such that $\mu(\omega)(B \Delta C) < \varepsilon$. Hence,

$$\begin{aligned} V_n^B(\omega) &\leq \sup_{A \in \mathcal{F}_n} |\mu(\omega)(A \cap B) - \mu(\omega)(A \cap C)| \\ &\quad + \sup_{A \in \mathcal{F}_n} |\mu(\omega)(A \cap C) - \mu(\omega)(A)\mu(\omega)(C)| \\ &\quad + \sup_{A \in \mathcal{F}_n} |\mu(\omega)(A)\mu(\omega)(C) - \mu(\omega)(A)\mu(\omega)(B)| \\ &\leq V_n^C(\omega) + 2\mu(\omega)(B \Delta C) < V_n^C(\omega) + 2\varepsilon \qquad \text{for all } n. \end{aligned}$$

Since $\omega \in A_0$ and $C \in \mathcal{F}_1$, it follows that

$$\limsup_n V_n^B(\omega) \leq 2\varepsilon + \limsup_n V_n^C(\omega) = 2\varepsilon \qquad \text{for all } B \in \mathcal{A}_1 \text{ and } \varepsilon > 0.$$

Therefore, $\mu(\omega)$ is 0–1 on $\mathcal{A}$. It remains to check condition (7). Fix $B \in \mathcal{A}_1$, take any version of $E(I_B \mid \mathcal{A}_n)$ and define $Z_n = E(I_B \mid \mathcal{A}_n) - \mu(B)$. Then, $|Z_n| \leq 2$ a.s. for all $n$, and the martingale convergence theorem yields $Z_n \to 0$ a.s. Further, for fixed $n \geq 1$ and $A \in \mathcal{F}_n$, one obtains

$$\begin{aligned} |E(I_A I_B \mid \mathcal{A}) &- E(I_A \mid \mathcal{A}) E(I_B \mid \mathcal{A})| \\ &= |E(I_A E(I_B \mid \mathcal{A}_n) \mid \mathcal{A}) - E(I_A E(I_B \mid \mathcal{A}) \mid \mathcal{A})| \\ &= |E(I_A Z_n \mid \mathcal{A})| \leq E(|Z_n| \mid \mathcal{A}) \qquad \text{a.s.} \end{aligned}$$

Since $\mathcal{F}_n$ is countable, it follows that

$$V_n^B = \sup_{A \in \mathcal{F}_n} |E(I_A I_B \mid \mathcal{A}) - E(I_A \mid \mathcal{A}) E(I_B \mid \mathcal{A})| \leq E(|Z_n| \mid \mathcal{A}) \to 0 \qquad \text{a.s.} \square$$

As noted in Section 2.2, a tail sub-$\sigma$-field includes its atoms so that (4) implies (3). Thus, by Theorem 15, condition (3) holds provided $\mathcal{A}$ is a tail sub-$\sigma$-field and $P$ admits a r.c.d. $\mu$ given $\mathcal{A}$, even if the other assumptions of Lemma 7 fail. [In fact, such assumptions grant something more than (3)].

**Acknowledgment.** We are indebted to Michael Wichura for various helpful suggestions and for having improved Theorems 3 and 8.

Dipartimento di Matematica  
Pura ed Applicata "G. Vitali"  
Universita' di Modena e Reggio-Emilia  
via Campi 213/B  
41100 Modena  
Italy  
E-mail: berti.patrizia@unimore.it

Dipartimento di Economia Politica  
e Metodi Quantitativi  
Universita' di Pavia  
via S. Felice 5  
27100 Pavia  
Italy  
E-mail: prigo@eco.unipv.it